\documentclass[12pt,thmsa]{article}
\usepackage{amsmath}
\usepackage{graphicx}
\usepackage{amsfonts}
\usepackage{amssymb}
\usepackage{mathrsfs}
\newcommand{\R}{\mathbb{R}}
\newcommand{\N}{\mathbb{N}}
\newcommand{\E}{\mathbb{E}}
\newcommand{\p}{\mathbb{P}}
\newcommand{\Z}{\mathbb{Z}}

\newtheorem{thm}{Theorem}
\newtheorem{theorem}{Theorem}

\newtheorem{lemma}{Lemma}
\newtheorem{corollary}{Corollary}

\newtheorem{remark}{Remark}

\newtheorem{Assumption}{Assumption}
\begin{document}

\begin{center}
{\Large \textbf{On estimation and prediction in a spatial semi-functional linear regression model}}

\bigskip

S\'ephane  BOUKA$^{a,\ast}$,  Kowir PAMBO BELLO$^{b,\dagger} $ and  Guy Martial  NKIET$^{a,\ddagger}$ 

\bigskip

$^a$Universit\'{e} des Sciences et Techniques de Masuku,  Franceville, Gabon.

\medskip

$^b$Ecole Normale Sup\'erieure,  Libreville, Gabon.

\bigskip

E-mail adresses : $^\ast$stephane.bouka@univ-masuku.org,    $^\dagger$kowir.pambo-bello@ens-libreville.org, $^\ddagger$guymartial.nkiet@univ-masuku.org.

\bigskip
\end{center}

\noindent\textbf{Abstract.}  We tackle estimation and prediction at non-visted sites in a spatial semi-functional linear regression model with derivatives that combines  a functional linear model with a nonparametric regression one. The parametric part is estimated by a method of moments and the other one by a local linear estimator. We establish the convergence rate of  the resulting estimators and predictor. A simulation study and an application to ozone pollution prediction at non-visted sites are proposed to illustrate our results.
\bigskip


\section{Introduction}

\noindent For two decades, there has been increased interest in statistical modelling of functional data. Among the methods allowing to deal with such data,   the functional  linear   regression, discussed in \cite{cardot99,cardot03},  constitutes an essential tool  and has a large number of applications in several fields  (see, e.g., \cite{ramsay97}). The corresponding model has seen some recent developments. Indeed, in order to enhance prediction power, \cite{mas_pumo09}  added  a component that   includes a derivative.  Furthermore, other works considered the so-called semi-functional partially linear regression model that combines a functional linear model with a nonparametric regression model. 
Estimation  based on this previous class of   models  has been   investigated (see \cite{lingetal19}, \cite{zhou_chen12}, \cite{zhouetal16}, and \cite{zhuetal20}), but only for purely non-spatial data. Despite the interest of processing spatial data, there are few works dealing with this type of data, compared to those of non-spatial data. However, \cite{giraldoetal11}  proposed a  methodology to make spatial linear predictions at non-data locations when the data are functions, which is applied to predict a curve of temperature; and \cite{bouka3}  considered spatial functional regression model with a derivative, which can be applied to predict the ozone pollution at non-visited sites. On the other hand, fixed design nonparametric regression  has also received a special attention in the field of spatial literature (see \cite{bouka}, \cite{mach_stoica10}, \cite{wang_wang}). Some results from this class of models can be applied to image analysis (see for instance \cite{mach_stoica10}). However, to the best of our knowledge, a small attention has been taken for the estimation in semi-functional linear regression model for spatially dependent observations (see for instance \cite{huetal20}, \cite{li_ying21} and \cite{benallouetal21}).

Consider $G:=L^2[0,1]$ and   the  Sobolev space  $H=\{x\in G, x^\prime\in G \}$, where   $x^\prime$ is  the first derivative of $x$; these are Hilbert spaces with inner products  $\left\langle .,.\right\rangle_{H}$ and $\left\langle .,.\right\rangle_{G}$  defined by:
\begin{eqnarray*}
\left\langle f, g\right\rangle_{G}=\int^{1}_{0}f(t)g(t)dt,\,\,\,\left\langle f, g\right\rangle_{H}=<f,g>_G+<f^\prime,g^\prime>_G.
\end{eqnarray*} 
We are interested with the model: 
\begin{eqnarray}\label{1.2}
Y_{\mathbf{i}}=\left\langle \phi, X_{\mathbf{i}}\right\rangle_{H}+\left\langle \gamma, X^{\prime}_{\mathbf{i}}\right\rangle_{G}+r\left(\frac{i_1}{n_1+1},\cdots,\frac{i_d}{n_d+1}\right)+\epsilon_{\mathbf{i}}\ ,
\end{eqnarray} 
where $\mathbf{i}=(i_1,\cdots,i_d)\in\mathcal{I}_{\mathbf{n}}:=\{\mathbf{i}=(i_1,\cdots,i_d)\in\Z^d: 1\leq i_k\leq n_k, k=1,\cdots,d\},\ d\ge 2$, $\phi$ and $\gamma$ are functions, $r(.)$ is a nonparametric spatial function defined on $[0,1]^d$ to be estimated by local linear smoothing, the errors $\epsilon_{\mathbf{i}}$ are spatially correlated, with a covariance function as assumed in Assumption {\ref{rega1}}. The random functions $X_{\mathbf{i}}$ and $X^{\prime}_{\mathbf{i}}$ are assumed to be centered and independent of the noise $\epsilon_{\mathbf{i}}$ which is also centered. The triplet $(X_{\mathbf{i}},Y_{\mathbf{i}},X^{\prime}_{\mathbf{i}})$ has the same probability distribution than the random vector $(X,Y,X^{\prime})$.  We are interested  in the prediction at a non-visited site obtained from the estimation of  the unknown parameter $(\phi, \gamma, r)$ defined in (\ref{1.2}). The special case with $\phi=0$ and $\gamma=0$ is considered in \cite{mach_stoica10},\cite{bouka}, whereas that with $r=0$ is studied in \cite{bouka3}.  
The article is organized as follows. Our estimation procedure is given in Section \ref{s2}. Section \ref{s3} is devoted to the assumptions and the main  results. In Section \ref{s4}, a  simulation study is presented whereas an application to ozone pollution forecasting at the non-visited sites is presented to Section \ref{s5}. A discussion of the results is made in Section \ref{Discussion}, while proofs of asymptotic results are postponed in Section  \ref{s6}.

\section{Estimation}\label{s2} 

\noindent Our estimation procedure is an association between the method of moments proposed in \cite{mas_pumo09} and the one based on local linear approximation used in \cite{bouka}. For estimating  the pair ($\phi$, $\gamma$), we adopt a method similar to that of \cite{mas_pumo09}. We first  multiply both members of   (\ref{1.2}) by $\left\langle X_{\mathbf{i}}, . \right\rangle_{H}$ and take the expectation of the obtained term.  Secondly, we reproduce this procedure  by multiplying by $\left\langle X^{\prime}_{\mathbf{i}}, . \right\rangle_{G}$. Since $X_{\mathbf{i}}$ and $X^{\prime}_{\mathbf{i}}$ are assumed to be centered, we finally obtain the following system:
\begin{eqnarray}\label{1.3}
\left\lbrace
\begin{array}{l}
\Delta=\Gamma \phi+\Gamma^{\prime}\gamma\\

\Delta^{\prime}=(\Gamma^{\prime } )^\ast\phi+\Gamma^{\prime\prime}\gamma
\end{array},
\right.
\end{eqnarray} 
where $\Gamma=\E(X\otimes_{H}X)$,   
$\Gamma^{\prime}=\E(X^{\prime}\otimes_{G}X)$, $\Gamma^{\prime\prime}=\E(X^{\prime}\otimes_{G}X^{\prime})$, $\Delta=\E(YX)$, $\Delta^{\prime}=\E(Y  X^{\prime})$, $(\Gamma^{\prime } )^\ast$ is the adjoint of $\Gamma^{\prime }$,  and where $\otimes_H$ (resp. $\otimes_G$) denotes the tensor product defined by $(u\otimes_H v)(h)=\left\langle u,h\right\rangle_H v$ (resp. $(u\otimes_G v)(h)=\left\langle u,h\right\rangle_G v$). The solution of the system (\ref{1.3}) is given by $\phi=S_{\phi}^{-1}\left[\Delta-\Gamma^{\prime}\Gamma^{\prime\prime-1}\Delta^{\prime}\right]$ and $  \gamma=S_{\gamma}^{-1}\left[\Delta^{\prime}-\Gamma^{\prime *}\Gamma^{-1}\Delta\right]$,  where $S_{\phi}=\Gamma-\Gamma^{\prime}\Gamma^{\prime\prime-1}\Gamma^{\prime*}$, $S_{\gamma}=\Gamma^{\prime\prime}-\Gamma^{\prime*}\Gamma^{-1}\Gamma^{\prime}$.

\noindent Considering the empirical versions of these operators and functions given by
\begin{eqnarray*}\label{emp1}
\Gamma_{\mathbf{n}}&=&\frac{1}{\widehat{\mathbf{n}}}\sum_{\mathbf{i}\in{\mathcal I}_{\mathbf{n}}}X_{\mathbf{i}}\otimes_{H}X_{\mathbf{i}},\     
\Gamma^{\prime}_{\mathbf{n}}=\frac{1}{\widehat{\mathbf{n}}}\sum_{\mathbf{i}\in{\mathcal I}_{\mathbf{n}}}X^{\prime}_{\mathbf{i}}\otimes_{G}X_{\mathbf{i}},\  
\Delta_{\mathbf{n}}=\frac{1}{\prod_{k=1}^d(n_k-1)}\sum_{\mathbf{i}\in{\mathcal I}^{1}_{\mathbf{n}}}Y_{\mathbf{i}}X_{\mathbf{i}},\\
\Gamma^{\prime *}_{\mathbf{n}}&=&\frac{1}{\widehat{\mathbf{n}}}\sum_{\mathbf{i}\in{\mathcal I}_{\mathbf{n}}}X_{\mathbf{i}}\otimes_{H}X^{\prime}_{\mathbf{i}},\     
\Gamma^{\prime\prime}_{\mathbf{n}}=\frac{1}{\widehat{\mathbf{n}}}\sum_{\mathbf{i}\in{\mathcal I}_{\mathbf{n}}}X^{\prime}_{\mathbf{i}}\otimes_{G}X^{\prime}_{\mathbf{i}},\  
\Delta^{\prime}_{\mathbf{n}}=\frac{1}{\prod_{k=1}^d(n_k-1)}\sum_{\mathbf{i}\in{\mathcal I}^{1}_{\mathbf{n}}}Y_{\mathbf{i}}X^{\prime}_{\mathbf{i}},
\end{eqnarray*}
where ${\mathcal I}^{1}_{\mathbf{n}}=\prod_{k=1}^d\{2,\cdots,n_k\}$, $\mathbf{n}=(n_1,\cdots,n_k)$, $\widehat{\mathbf{n}}=\prod_{k=1}^dn_k$. We write $\mathbf{n}\rightarrow +\infty$ if $\min_{1\le k\le d}(n_k)\longrightarrow +\infty$ and we assume that $\dfrac{n_j}{n_k}\le C$ for $1\le j,\ k\le d$ and $0<C<+\infty$. Invertible empirical operators are obtained by regularization method: 
$$ \widetilde{\Gamma}^{-1}_{\mathbf{n}}=(\Gamma_{\mathbf{n}}+w_{\mathbf{n}}I)^{-1},\, \widetilde{\Gamma}^{\prime\prime -1}_{\mathbf{n}}=(\Gamma^{\prime\prime}_{\mathbf{n}}+w_{\mathbf{n}}I)^{-1},\,  $$ 
where $w_{\mathbf{n}}$  is a  sequence, from $\N^d$ to $\R$, converging to $0$ as  $\mathbf{n}\rightarrow +\infty$  and $I$ denotes the identity operator. We put $$S_{\mathbf{n},\phi}=\Gamma_{\mathbf{n}}-\Gamma^{\prime}_{\mathbf{n}}(\widetilde{\Gamma}^{\prime\prime -1}_{\mathbf{n}})\Gamma^{\prime*}_{\mathbf{n}},\ 
S_{\mathbf{n},\gamma}=\Gamma^{\prime\prime}_{\mathbf{n}}-\Gamma^{\prime*}_{\mathbf{n}}(\widetilde{\Gamma}^{-1}_{\mathbf{n}})\Gamma^{\prime}_{\mathbf{n}},$$
$$\ \ u_{\mathbf{n},\phi}=\Delta_{\mathbf{n}}-\Gamma^{\prime}_{\mathbf{n}}(\widetilde{\Gamma}^{\prime\prime -1}_{\mathbf{n}})\Delta^{\prime}_{\mathbf{n}},\ 
 u_{\mathbf{n},\gamma}=\Delta^{\prime}_{\mathbf{n}}-\Gamma^{\prime*}_{\mathbf{n}}(\widetilde{\Gamma}^{-1}_{\mathbf{n}})\Delta_{\mathbf{n}},$$ 
and we  estimate $(\phi,\gamma)$ by the pair $(\widehat{\phi}_{\mathbf{n}},\widehat{\gamma}_{\mathbf{n}})$  given, as in \cite{bouka3}, by:
\begin{equation*}\label{2.2}
\widehat{\phi}_{\mathbf{n}}=(S_{\mathbf{n},\phi}+\psi_{\mathbf{n}}I)^{-1}u_{\mathbf{n},\phi},\ \ 
\widehat{\gamma}_{\mathbf{n}}=(S_{\mathbf{n},\gamma}+\psi_{\mathbf{n}}I)^{-1}u_{\mathbf{n},\gamma},
\end{equation*}
where $\psi_{\mathbf{n}}$  is a  sequence, from $\N^d$ to $\R$, converging to $0$ as  $\mathbf{n}\rightarrow +\infty$.
 Secondly,  for estimating  the nonparametric regression function, we rewrite  (\ref{1.2}) as  
\begin{eqnarray}\label{1.4}
T_{\mathbf{i}}=r\left(\frac{\mathbf{i}}{\mathbf{n+1}}\right)+\xi_{\mathbf{i}},
\end{eqnarray}
where $\dfrac{\mathbf{i}}{\mathbf{n+1}}=\left(\dfrac{i_1}{n_1+1},\cdots,\dfrac{i_d}{n_d+1}\right)$,  $T_{\mathbf{i}}=Y_{\mathbf{i}}-\left\langle \widehat{\phi}_{\mathbf{n}}, X_{\mathbf{i}}\right\rangle_{H}-\left\langle \widehat{\gamma}_{\mathbf{n}}, X^{\prime}_{\mathbf{i}}\right\rangle_{G}$ and $\xi_{\mathbf{i}}=\epsilon_{\mathbf{i}}+\left\langle \phi-\widehat{\phi}_{\mathbf{n}}, X_{\mathbf{i}}\right\rangle_{H}+\left\langle \gamma-\widehat{\gamma}_{\mathbf{n}}, X^{\prime}_{\mathbf{i}}\right\rangle_{G}$, and we   locally approximate $r$  by a linear regression function by using Taylor expansion in the neighbourhood of $\mathbf{s}_{0}\in[0,1]^{d}$ i.e. $r(\frac{\mathbf{i}}{\mathbf{n+1}})\approx\beta_{0}+\beta^{T}_{1}(\frac{\mathbf{i}}{\mathbf{n+1}}-\mathbf{s}_{0})$, 
where $\beta_{0}=r(\mathbf{s}_{0})$ and $\beta_{1}=\nabla r(\mathbf{s}_{0})$ is the gradient vector of $r$ at $\mathbf{s}_{0}$. From the model defined in (\ref{1.4}), an estimator of $\beta_{0}$ is given by the solution of the following least squares minimization problem:
\begin{eqnarray*}
\min_{\beta_{0},\beta_{1}}\sum_{\mathbf{i}\in{\mathcal I}_{\mathbf{n}}}\left\{T_{\mathbf{i}}-\beta_{0}-\beta^{T}_{1}\left(\frac{\mathbf{i}}{\mathbf{n+1}}-\mathbf{s}_{0}\right)\right\}^2\frac{1}{h^{d}}K\left(\frac{\frac{\mathbf{i}}{\mathbf{n+1}}-\mathbf{s}_{0}}{h}\right),
\end{eqnarray*}
where $K:\R^{d}\longrightarrow\R_{+}$ is a kernel function and $h$ is a bandwidth, and $u^T$ is the transpose of $u$. Then, a local linear estimator $\widehat{r}$ of $r$ on $\mathbf{s}_0\in[0,1]^d$ is given,  as in \cite{bouka}, by:
\begin{eqnarray*}
\widehat{r}(\mathbf{s}_{0})=(1, \mathbf{0}^{T})\left(\frac{1}{\widehat{\mathbf{n}}}{\mathcal X}^{T}W_{0}{\mathcal X}\right)^{-1}\left(\frac{1}{\widehat{\mathbf{n}}}{\mathcal X}^{T}W_{0}{\mathcal Y}\right)=\mathcal{S}^T_{\mathbf{s}_{0}}{\mathcal Y},
\end{eqnarray*}
where $(1,\mathbf{0}^T)\in\R^{d+1}$, $\mathcal{X}$, $\mathcal{Y}$ and $W_0$ are, respectively, the $\widehat{\mathbf{n}}\times (d+1)$, $\widehat{\mathbf{n}}\times 1$ and $\widehat{\mathbf{n}}\times \widehat{\mathbf{n}}$ matrices given by:
\[
\mathcal{X}=
\begin{pmatrix}
	1&\left(\dfrac{\dfrac{\mathbf{1}}{\mathbf{n+1}}-\mathbf{s}_{0}}{h}\right)^{T}\\
	\vdots&\vdots&\\
	1&\left(\dfrac{\dfrac{\mathbf{n}}{\mathbf{n+1}}-\mathbf{s}_{0}}{h}\right)^{T}\\
\end{pmatrix} ,\,\,\,\,\,\, \mathcal{Y}=\left(T_{\mathbf{1}},\cdots,T_{\mathbf{n}}\right)^{T} ,
\]
where $\mathbf{1}=(1,\cdots,1)\in\R^d$, $\mathbf{n}=(n_1,\cdots,n_d)$ and
 $$  W_{0}=\textrm{diag}\left\{\frac{1}{h^{d}}K\left(\dfrac{\dfrac{\mathbf{1}}{\mathbf{n+1}}-\mathbf{s}_{0}}{h}\right),\cdots,\frac{1}{h^{d}}K\left(\dfrac{\dfrac{\mathbf{n}}{\mathbf{n+1}}-\mathbf{s}_{0}}{h}\right)\right\}.$$



\section{Assumptions and results}\label{s3}

\noindent We first make assumptions needed for establishing our results. For $\beta>0$, $L>0$ and $M>0$ we denote by $\mathcal{H}(\beta, L)$ (resp.   $\mathcal{G}(L)$) the  H\"older class (resp. the class)  of functions $f$ satisfying $|f(x)-f(y)|\le L\left\|x-y\right\|^\beta_{\infty}$ (resp. $|f(x)-f(y)|\le L\left\|x-y\right\|_{\infty}$), and we consider the class:

 \begin{eqnarray*}\label{def}
 \Sigma(\beta, L, M)=\left\lbrace
 \begin{array}{l}
 \left\{f\in \mathcal{H}(\beta, L):\left\|f\right\|_{\infty}\leq M\right\}\  \text{if}\  \frac{d}{4}<\beta\leq 1\\
  \\
 \left\{f\in\mathcal{G}(L):\left\|f\right\|_{\infty}\leq M, \left\|\nabla f\right\|_{\infty}\leq M\right\}\  \text{if}\  \beta>1
 \end{array}.
 \right. 
 \end{eqnarray*}
  
 \begin{Assumption}\label{as0} 
 The function $r$ belongs to  $\Sigma(\beta, L, M)$ with $\beta>d/4$, $L>0$ and $M>0$.
  \end{Assumption}

  \begin{Assumption}\label{as1}
  $Ker(\Gamma)=Ker(\Gamma^{\prime\prime})=\{0\}$,
  where  $Ker (A)=\{x \,: A x=0\}$. 
  \end{Assumption}
  
  \begin{Assumption}\label{as2}
  $(\phi,\gamma)\notin\{(f,g)\in H\times G: f+D^{*}g=0\}$ where $D^{*}$ is the adjoint of the ordinary differential operator $D$. 
  \end{Assumption}

  \begin{Assumption}\label{as3}
  $\left\|\Gamma^{-1/2}\phi\right\|_{H}<+\infty$,\ \ 
  $\left\|(\Gamma^{\prime\prime})^{-1/2}\gamma\right\|_{G}<+\infty$, where $\Vert\cdot\Vert_H$ and $\Vert\cdot\Vert_G$ are the norms induced by $<\cdot ,\cdot>_H$ and  $<\cdot ,\cdot>_G$  respectively.
  \end{Assumption}

  \begin{Assumption}\label{as3.1}
  The process $\{Z_{\mathbf{i}}=(X_{\mathbf{i}}, Y_{\mathbf{i}}, X^{\prime}_{\mathbf{i}}),\; \mathbf{i}\in\mathbb{Z}^d\}$ is $\alpha$-mixing dependent. That is\\ $\lim_{m\rightarrow +\infty}\alpha_{1,\infty}(m)=0$, where 
  \begin{eqnarray}\label{ar2}
  \alpha_{1,\infty}(m)=\sup_{\{\mathbf{i}\}, E\subset\Z^{d}, \rho(E,\{\mathbf{i}\})\geq m}\{\sup_{A\in\sigma(Z_{\mathbf{i}}),B\in\sigma(Z_{\mathbf{j}};\mathbf{j}\in{E})}\{|\p(A\cap B)-\p(A)\p(B)|\}\},
  \end{eqnarray}
   $\rho$ is the distance defined for any subsets $E_{1}$ and $E_{2}$ of $\Z^{d}$, by $\rho(E_{1},E_{2})=\min\{\|\mathbf{i}-\mathbf{j}\|, \mathbf{i}\in E_{1}, \mathbf{j}\in E_{2}\}$ with $\|\mathbf{i}-\mathbf{j}\|=\max_{1\leq s\leq d}|i_{s}-j_{s}|$.
  \end{Assumption}
  \begin{Assumption}\label{as5}
  $\left\|X_{\mathbf{i}} \right\|_{H}\leq C$ a.s. where $C$ is some positive constant.
  \end{Assumption}
  
 \begin{Assumption}\label{rega1}
 $Cov(\epsilon_\mathbf{i},\epsilon_\mathbf{j})=\sigma^2\exp(-a\|\mathbf{i}-\mathbf{j}\|)$ for all  $\mathbf{i}, \mathbf{j}\in{\mathcal I}_{\mathbf{n}}$, where $a$ and $\sigma^2$ are known positive constants.
 \end{Assumption}
 \begin{Assumption}\label{rega2}
 The kernel function $K(.)$ is symmetric, Lipschitz, continuous and bounded. The support of $K(.)$ is $[-1,1]^{d}$, $\int K(\mathbf{u})d{\mathbf{u}}=1$, $\int {\mathbf{u}}K(\mathbf{u})d{\mathbf{u}}=\mathbf{0}$, $\int {\mathbf{u}}{\mathbf{u}}^{\tau}K(\mathbf{u})d{\mathbf{u}}=\nu_{2}I_d$, where $I_d$ is the $d\times d$ identity matrix  and $\nu_{2}\ne0$.
 \end{Assumption}
 
 \begin{remark} \label{cor1}
 \noindent Assumptions \ref{as1}--\ref{as5} are technical conditions to ensure consistency of $\widehat{\phi}_{\mathbf{n}}$ and $\widehat{\gamma}_{\mathbf{n}}$ (see \cite{bouka3}). However, Assumption \ref{as5} can be replaced by $\E\left(\left\|X\right\|^8_H\right)<C$ (see \cite{mas_pumo09}), but this assumption on finite moment would lead us to longer and more intricate
 methods of proof.  Assumptions \ref{rega1} et \ref{rega2} are conditions needed to establish   consistency of the estimator $\widehat{r}$ of $r $; they have also been used in \cite{fransisco} and in \cite{bouka}. However, spatial covariance models such that Matern's spatial covariance model or Gaussian covariance model can be also used. 
  \end{remark} 
  
\bigskip

 \noindent Let $\Lambda_{\mathbf{n}}$ (resp. $\Lambda$)  be one of the following:  $\Gamma_{\mathbf{n}}$, $\Gamma^{\prime}_{\mathbf{n}}$, 
    $\Gamma^{\prime*}_{\mathbf{n}}$,  $\Gamma^{\prime\prime}_{\mathbf{n}}$, $\Delta_{\mathbf{n}}$ and $\Delta^{\prime}_{\mathbf{n}}$ (resp. $\Gamma,\Gamma^{\prime},\Gamma^{\prime*},\Gamma^{\prime\prime}$, $\Delta$ and $\Delta^{\prime}$). Rates of convergence of $\Lambda_{\mathbf{n}}$ with respect to the norms $\|.\|_{\infty}$ and $\|.\|_{L^{2}({\mathcal HS})}$ are given in Theorem \ref{th1}, whereas those of $\widehat{\phi}_{\mathbf{n}}$ and $\widehat{\gamma}_{\mathbf{n}}$ are given in Corollary \ref{coro1}.

\bigskip
  
  \begin{thm}\label{th1}
   Let $(v_{j})_{j\geq1}$ (resp. $(v_{j}^{\star})_{j\geq1}$) be a complete orthonormal system in $H$ (resp. $G$) and $(\lambda_{j})_{j\geq1}$ be a characteristic numbers sequence of $\Lambda$ (i.e. the square root of the eigenvalues of $\Lambda^{*}\Lambda$)  with $\lambda_{j}=O(u^{j})$, $0<u<1$, $j\geq1$. Under assumptions $\ref{as1}-\ref{as5}$ with $\alpha_{1,\infty}(t)=O(t^{-\theta})$, $\theta>2d$, we have, for all $\tau>0$: 
  \begin{eqnarray}
  &&\p\left(\left\|\Lambda_{\mathbf{n}}-\Lambda\right\|_{\infty}>\tau\right)=O\left(\dfrac{\log \widehat{\mathbf{n}}}{\widehat{\mathbf{n}}}\right);\label{3.2}\\
  &&\left\|\Lambda_{\mathbf{n}}-\Lambda\right\|_{L^{2}({\mathcal HS})}=O\left(\dfrac{\log \widehat{\mathbf{n}}}{\sqrt{\widehat{\mathbf{n}}}}\right) \label{3.3}.
  \end{eqnarray}
  \end{thm}

\bigskip
  
\noindent From Theorem \ref{th1} and arguing as in \cite{mas_pumo09}, we derive the following corollary.
  \begin{corollary}\label{coro1}
   Under assumptions of Theorem \ref{th1}, we have: 
    $$\E\left(\left\|\phi-\widehat{\phi}_{\mathbf{n}}\right\|^2_{\Gamma}\right)=O\left(\frac{\psi^{2}_{\mathbf{n}}}{w^{2}_{\mathbf{n}}}\right)+O\left(\frac{(\log \widehat{\mathbf{n}})^{2}}{w^{2}_{\mathbf{n}}\psi^{2}_{\mathbf{n}}\widehat{\mathbf{n}}}\right);$$
    $$  \E\left(\left\|\gamma-\widehat{\gamma}_{\mathbf{n}}\right\|^2_{\Gamma^{\prime\prime}}\right)=O\left(\frac{\psi^{2}_{\mathbf{n}}}{w^{2}_{\mathbf{n}}}\right)+O\left(\frac{(\log \widehat{\mathbf{n}})^{2}}{w^{2}_{\mathbf{n}}\psi^{2}_{\mathbf{n}}\widehat{\mathbf{n}}}\right),$$  
    where $\|.\|_{\Gamma}:=\|\Gamma^{1/2}(.)\|_{H}$ and $\|.\|_{\Gamma^{\prime\prime}}:=\|\Gamma^{\prime\prime1/2}(.)\|_{G}$ are two semi-norms.
  \end{corollary}

\bigskip
  
\noindent  It remains to obtain bounds of the estimator $\widehat{r}$ of $r$.  The following theorem gives bounds of the bias and of the variance.

\bigskip
    
    \begin{theorem}\label{regl2}
    Assume that Assumptions \ref{as0}--\ref{rega2} are satisfied with $\alpha_{1,\infty}(t)=O(t^{-\theta})$, $\theta> 2(d+1)$. Moreover, assume that the eigenvalues $\lambda_j$ of the operator $\Gamma$ are such that $\lambda_j=O(u^j)$ with $0<u<1$, $j\ge1$, that $h\to 0$ and $\min_{k=1,\cdots,d}\{n_k\}h\to +\infty$ as $\mathbf{n}\to+\infty$. Then: 
   \item[(i)]$$   \sup_{\mathbf{s}_{0}\in[0,1]^d}\left[\E(\widehat{r}(\mathbf{s}_{0}))-r(\mathbf{s}_{0})\right]^2=O\left(h^4\right)+O\left(\frac{(\log \widehat{\mathbf{n}})^2}{w_{\mathbf{n}}^2\psi_{\mathbf{n}}^2\widehat{\mathbf{n}}}\right)+O\left(\psi_{\mathbf{n}}^2\right);$$
   \item[(ii)]$$ \sup_{\mathbf{s}_{0}\in[0,1]^d}Var(\widehat{r}(\mathbf{s}_{0}))=O\left(\frac{\log \widehat{\mathbf{n}}}{\widehat{\mathbf{n}}h^{d}}\right).$$ 
    \end{theorem}

\bigskip
   
   \begin{remark}\label{cor2}
   An immediate consequence of Theorem \ref{regl2} is for all $\mathbf{s}_{0}\in[0,1]^d$, $$\left|\widehat{r}(\mathbf{s}_{0})-r(\mathbf{s}_{0})\right|=O_p\left(\dfrac{(\log \widehat{\mathbf{n}})^{2}}{w^{2}_{\mathbf{n}}\psi^{2}_{\mathbf{n}}\widehat{\mathbf{n}}}+h^4+\psi_{\mathbf{n}}^2+\dfrac{\log \widehat{\mathbf{n}}}{\widehat{\mathbf{n}}h^d}\right)$$
    and so the optimal bandwidth $h$ is controlled by the trade-off between the variance and the square of the bias. If $\psi_{\mathbf{n}}\propto \dfrac{(\log \widehat{\mathbf{n}})^{1/2}}{\widehat{\mathbf{n}}^{2/(4+d)}}$, $h=\dfrac{1}{\widehat{\mathbf{n}}^{1/(4+d)}}$ with $d\ge 2$, and $\dfrac{1}{w_{\mathbf{n}}^2}\propto\log\widehat{\mathbf{n}}$, then for all $\mathbf{s}_{0}\in[0,1]^d$, we have \[\left|\widehat{r}(\mathbf{s}_{0})-r(\mathbf{s}_{0})\right|=O_p\left(\max\left(\frac{\left(\log \widehat{\mathbf{n}}\right)^2}{\widehat{\mathbf{n}}^{d/(4+d)}},\frac{\log \widehat{\mathbf{n}}}{\widehat{\mathbf{n}}^{4/(4+d)}}\right)\right)\] which is, with $d\ge4$, the optimal convergence rate in spatial nonparametric regression setting when data are $\alpha$-mixing dependent and the correlation of errors is long-range, and for $d=3$, it is better than $O_p\left(\left(\frac{\log \widehat{\mathbf{n}}}{\widehat{\mathbf{n}}^{4/(4+d)}}\right)^{1/2}\right)$ given in \cite{carbonetal07}, and is quite close from the one of \cite{carbonetal07} with $d=2$. Besides, when the correlation of errors is short-range, the optimal convergence rate from $\widehat{r}$ to $r$ is $O_p\left(\dfrac{1}{\widehat{\mathbf{n}}^{4/(4+d)}}\right)$  \cite[p. 36]{liu01}.
    \end{remark}

\bigskip
    
\noindent In addition, from Remark \ref{cor1}, Theorem \ref{regl2}, Lemma \ref{l1} (see Section \ref{s5}) and Theorem \ref{th1} together with the same arguments than in \cite{mas_pumo09}, we deduce the following Corollary.

\bigskip
      
    \begin{corollary}
    Under Assumptions of Theorem \ref{regl2}, we have for each $\mathbf{i}_0\in{\mathcal I}_{\mathbf{n+1}}\setminus{\mathcal I}_{\mathbf{n}}$ that:$$\E\left[\left(\widehat{Y}_{\mathbf{i}_0}-Y_{\mathbf{i}_0}^*\right)^2\right]=O\left(\frac{\psi^{2}_{\mathbf{n}}}{w^{2}_{\mathbf{n}}}+\frac{(\log \widehat{\mathbf{n}})^{2}}{w^{2}_{\mathbf{n}}\psi^{2}_{\mathbf{n}}\widehat{\mathbf{n}}}+h^4+\frac{\log \widehat{\mathbf{n}}}{\widehat{\mathbf{n}}h^d}\right),$$ where $\widehat{Y}_{\mathbf{i}_0}=\left\langle \widehat{\phi}_{\mathbf{n}}, X_{\mathbf{i}_0}\right\rangle_H+\left\langle \widehat{\gamma}_{\mathbf{n}}, X^{\prime}_{\mathbf{i}_0}\right\rangle_G+\widehat{r}\left(\dfrac{\mathbf{i}_0}{\mathbf{n+1}}\right)$ and $Y_{\mathbf{i}_0}^*=\left\langle \phi, X_{\mathbf{i}_0}\right\rangle_H+\left\langle \gamma, X^{\prime}_{\mathbf{i}_0}\right\rangle_G+r\left(\dfrac{\mathbf{i}_0}{\mathbf{n+1}}\right)$.
    \end{corollary} 

\bigskip
   
   \begin{remark}\label{rem3}
   The methodology proposed in this paper can also be applied when the design of the non-parametric regression function is random. That is $$Y_{\mathbf{i}}=\int_{0}^{1}\phi(t)X_{\mathbf{i}}(t)dt+\int_{0}^{1}\gamma(t)X_{\mathbf{i}}^{\prime}(t)dt+r(Z_{\mathbf{i}})+\epsilon_{\mathbf{i}}\ ,\, \mathbf{i}\in{\mathcal I}_{\mathbf{n}},$$
   where $Z_{\mathbf{i}}$ is a random vector independent of $X_{\mathbf{i}}$ and $X_{\mathbf{i}}^{\prime}$. An estimator $\widehat{r}$ of $r$ is obtained by considering, in our estimation procedure (see Section {\ref{s2}}), the product kernels (see \cite[chapter 3]{camille}) $K_1\left(\dfrac{x_{\mathbf{s}_0}-Z_{\mathbf{i}}}{h_1}\right)K_2\left(\dfrac{\mathbf{s}_0-\mathbf{i}}{(\mathbf{n+1})h_2}\right)$ instead of $K\left(\dfrac{\dfrac{\mathbf{i}}{\mathbf{n+1}}-\mathbf{s}_0}{h}\right)$. Besides, if the correlation of errors is of short-range with  $Cov(\epsilon_{\mathbf{i}},\epsilon_{\mathbf{j}})=\sigma^2\exp(-a\widehat{\mathbf{n}}\|\mathbf{i}-\mathbf{j}\|)$, then refer to \cite{fransisco} or \cite{liu01} for the convergence from $\widehat{r}$ to $r$. Moreover, if $\gamma=0$, we obtain a model which extends to spatial case that considered in \cite{zhou_chen12} with independent data.
   \end{remark}


\section{A simulation study}\label{s4}

\noindent In this section, we present a simulation study made in order to appreciate the finite sample performance of our proposal. We consider an equivalent
 definition to the model (\ref{1.2}) given by
 \begin{eqnarray}\label{simu1}
 Y_{\mathbf{i}}=\int_{0}^{1}\phi(t)X_{\mathbf{i}}(t)dt+\int_{0}^{1}\gamma(t)X_{\mathbf{i}}^{\prime}(t)dt+r\left(\frac{\mathbf{i}}{n+1}\right)+\epsilon_{\mathbf{i}},
 \end{eqnarray}
 where $d=2$, $n_1=n_2=n$ and $\mathbf{i}\in\{1,\cdots,n\}^2$. By using the lexicographic order in $\mathbb{Z}^2$, we simulated a sample  $\{(X_{\mathbf{i}_\ell},Y_{\mathbf{i}_\ell})\}_{1\leq  \ell\leq n^2}$ such that:
 \[
 X_{\mathbf{i}_\ell}(t)=\sum_{k=1}^{15}\Lambda_{\mathbf{i}_\ell,k}\,F_k(t) 
 \]
 where $F_1,\cdots,F_{15}$ are the  $15$-th  first elements of the Fourier basis,  the random vector $(\Lambda_{\mathbf{i}_1,k},\cdots,\Lambda_{\mathbf{i}_{n^2},k})^T$ is  obtained from a multivariate truncated normal distribution with zero mean, $n^2\times n^2$ covariance matrix $\Sigma^1$ with general term  $\Sigma^1_{ij}=\exp(-a\|\mathbf{i}_i-\mathbf{i}_j\|_2)$, where $a= 0.1,\ 1,\ 3,\ 200$, and with lower truncation limit $(0,\cdots,0)\in\R^{n^2}$ and upper truncation limit $(1,\cdots,1)\in\R^{n^2}$. When $a=200$, there is approximately no spatial autocorrelation in the process. The process is said strongly correlated when $a=0.1,\ 1$ and weakly correlated when $a=3$. The process $Y_{\mathbf{i}_\ell}$ is obtained from the model (\ref{simu1}), in which $X_{\mathbf{i}_\ell}^{\prime}$ is computed by the function "fdata.deriv" of the R fda package,   $\phi(t)=[\sin(2\pi t^3)]^3$, $\gamma(t)=(0.6-t)^2$, $t\in[0,1]$, $r(\mathbf{x})=\exp(-\|\mathbf{x}\|_{\infty})$, $\mathbf{x}\in[0,1]^2$, integrals are approximated by the rectangular method applied at the $366$ equispaced points of the interval $[0,1]$, $(\epsilon_{\mathbf{i}_1},\cdots,\epsilon_{\mathbf{i}_{n^2}})^T$ is a   random vector having a normal distribution ${\mathcal N}(0,\Sigma^2)$ where $\Sigma^2$ is a $n^2\times n^2$ covariance matrix with general term  $\Sigma^2_{ij}=0.01\Sigma^1_{ij}$ for $i\ne j$ and $\Sigma^2_{ii}=\Sigma^1_{ii}$. The estimate ($\widehat{\phi}_{\mathbf{n}}$, $\widehat{\gamma}_{\mathbf{n}}$) of the pair ($\phi$,$\gamma$) depends of the regularization sequences $\psi$ and $w$. These sequences are obtained from cross validation based on the mean standard error of prediction :
 $$CVMSEP(\psi,w)=\dfrac{1}{n^2}\sum^{n^2}_{\ell=1}\left(Y_{\mathbf{i}_{\ell}}-\widetilde{Y}^{(\ell)}_{\mathbf{i}_{\ell}}(\psi,w)\right)^2,$$
  where  $\widetilde{Y}^{(-\ell)}_{\mathbf{i}_{\ell}}(\psi,w)=\left\langle \widehat{\phi}_{\mathbf{n}}, X_{\mathbf{i}_{\ell}}\right\rangle_G+\left\langle \widehat{\gamma}_{\mathbf{n}}, X^{\prime}_{\mathbf{i}_{\ell}}\right\rangle_G$ with $\widehat{\phi}_{\mathbf{n}}$ and $\widehat{\gamma}_{\mathbf{n}}$  computed with the $\ell-$th part of the removed data. The estimate $\widehat{r}$ of $r$ depends on the bandwidth $h$ which is selected by minimizing the following generalized cross validation (GCV) function  \cite{fransisco}:$$GCV (h)=\dfrac{1}{n^2}\sum_{\ell=1}^{n^2}\left(\dfrac{T_{\mathbf{i}_{\ell}}-\widehat{r}\left(\frac{\mathbf{i}_{\ell}}{n};h\right)}{1-\frac{1}{n^2}tr\left(\mathbf{SC}\right)}\right)^2,$$
 where $T_{\mathbf{i}_{\ell}}=Y_{\mathbf{i}_{\ell}}-\left\langle \widehat{\phi}_{\mathbf{n}}, X_{\mathbf{i}_{\ell}}\right\rangle_G-\left\langle \widehat{\gamma}_{\mathbf{n}}, X^{\prime}_{\mathbf{i}_{\ell}}\right\rangle_G$ with $\widehat{\phi}_{\mathbf{n}}$ and $\widehat{\gamma}_{\mathbf{n}}$  computed from the optimal regularization parameters $\psi_{opt}$ and $w_{opt}$, $\widehat{r}\left(\frac{\mathbf{i}_{\ell}}{n};h\right)$ is  computed from Epanechnikov kernel defined by $K(x)=\frac{2}{\pi}\max\left\{(1-\|x\|_2^2),0\right\}$, $\mathbf{S}$ is the $n^2\times n^2$ matrix whose $\ell$th row is equal to $\mathcal{S}^T_{\mathbf{i}_{\ell}/n}$ and $\mathbf{C}$ is the correlation matrix of the observations. 
 
 We assess performance of our method through calculation of the Mean Squared Error ($MSE_1$ and $MSE_2$), based on $100$ replications with $n=5, 10$ and $a=0.1,\ 1,\ 3,\ 200$, and defined by:
   $$MSE_1=\frac{1}{n^2}\sum_{j=1}^{n}\sum_{i=1}^{n}\left[r\left(\frac{i}{n+1},\frac{j}{n+1}\right)-\widehat{r}\left(\frac{i}{n+1},\frac{j}{n+1}\right)\right]^2,$$
   \begin{eqnarray*}
    MSE_2&=&\frac{1}{n^2}\sum_{j=1}^{n}\sum_{i=1}^{n}\left[\left\langle \phi-\widehat{\phi}_{\mathbf{n}}, X_{(i,j)}\right\rangle_G+ \left\langle \gamma-\widehat{\gamma}_{\mathbf{n}}, X^{\prime}_{(i,j)}\right\rangle_G\right.\\
    &&\left.+ r\left(\frac{i}{n+1},\frac{j}{n+1}\right)-\widehat{r}\left(\frac{i}{n+1},\frac{j}{n+1}\right)\right]^2.
    \end{eqnarray*}  
    We denote by $m$ the mean and by $sd$ the standard deviation.
  The results are postponed in   Table \ref{tab:1}.
    
 \begin{table}[ht!]
  \begin{center}
       \begin{tabular}{cccccccccc}
       \hline
       \multicolumn{1}{c}{}&\multicolumn{1}{c}{}&\multicolumn{2}{c}{$a=0.1$}&\multicolumn{2}{c}{$a=1$}&\multicolumn{2}{c}{$a=3$}&\multicolumn{2}{c}{$a=200$}\\
       \cline{3-10}
        $n^2$& error criterion & $m$& $sd$& $m$& $sd$& $m$& $sd$& $m$& $sd$ \\\hline
       25&$MSE_1$ & 0.51& 0.58& 0.29& 0.17& 0.16& 0.10& 0.20& 0.09\\
       &$MSE_2$ & 0.98& 1.24& 0.31& 0.18& 0.18& 0.11& 0.21& 0.11\\
       & & & & & & & & &\\
       100&$MSE_1$ & 0.38& 0.40& 0.12& 0.08& 0.07& 0.04& 0.05& 0.03\\
       &$MSE_2$ & 0.62& 0.73& 0.15& 0.10& 0.06& 0.03& 0.06& 0.03\\\hline    
       \end{tabular}
       \caption{Mean ($m$) and standard deviation ($sd$) of $MSE_1$ and $MSE_2$, based on 100 replications; $n^2=25,\ 100$ and  $a =0.1,\ 1,\ 3,\ 200$.}
         \label{tab:1} 
       \end{center}
    \end{table}
    
\noindent In Table \ref{tab:1}, the two kinds of error criterion ($MSE_1$, $MSE_2$) have a decreasing general tendency  as the sample size increases. Thus our estimation procedure well fits to spatial semi-functional linear regression model with derivatives. Also, for each fixed $n$, values of each error criterion for weakly correlated processes ($a=3$) are similar to those approximatively non-correlated processes ($a=200$), whereas those of strongly correlated processes ($a=0.1,1$) decrease as $a$ increases, so showing the interest to consider spatially dependent observations in this study. Besides, for each fixed $n$, $MSE_1$ and $MSE_2$ are similar for values of $a\geq 1$. This means that the presence of functional data in the model (\ref{1.2}) does not modify the  convergence rate of the estimated nonparametric regression function as stated in \cite[Remark 2]{zhou_chen12}.


\section{Application to ozone pollution forecasting at the non-visited sites} \label{s5}

\noindent In this section, our methodology is applied to predict the level of ozone pollution at non-visited sites of California state. For that, we use the available data on internet site   https://www.epa.gov/outdoor-air-quality-data. The explicative functional variables $$\{X_{s_i}(t),\ t=1,\cdots,100,\ s_i=(Latitude, Longitude)_i,\  i=1,\cdots,51\}$$ correspond to the measurements of ozone concentration measured the $p=100$ firsts days,  from  January 1st, 2021 to April 12th, 2021 on each of $n^2=51$ sites.
 The response variables $$\{Y_{s_i},\ s_i=(Latitude, Longitude)_i,\ i=1,\cdots,35\}$$ correspond to the measurements of ozone concentration measured on April 13th, 2021 on each of $35$ firsts stations. $\{Y_{s_i},\ s_i=(Latitude, Longitude)_i,\ i=1,\cdots,35\}$ and $\{X_{s_i}(t),\ t=1,\cdots,100,\ s_i=(Latitude, Longitude)_i,\  i=1,\cdots,35\}$ are related through the spatial semi-functional linear regression model with derivatives (SSFLRD) defined by 
 \begin{eqnarray*}
 Y_{s_{i}}&=&\int_{0}^1\phi(t)X_{s_i}(t)dt+\int_{0}^1\gamma(t)X^{\prime}_{s_i}(t)dt\\
 &&+r\left(\frac{Latitude}{\max_{j=1,\cdots,51}(Latitude[j])},\frac{Longitude}{\max_{j=1,\cdots,51}(Longitude[j])}\right)+\epsilon_{s_i}.
 \end{eqnarray*}   
 For evaluating the performances of our method, we compare prediction error obtained from SSFLRD model with the spatial functional linear regression model with derivatives (SFLRD) studied in \cite{bouka3} and defined by
 $$Y_{s_{i}}=\int_{0}^1\phi(t)X_{s_i}(t)dt+\int_{0}^1\gamma(t)X^{\prime}_{s_i}(t)dt+\epsilon_{s_i}\ ,$$   
  where  $X^{\prime}_{\mathbf{i}_{\ell}}$ standing for the first derivative of $X_{\mathbf{i}_{\ell}}$ is computed from the function "fdata.deriv" of the $R$ fda package. So, we predict from both methodologies  $$\{Y_{s_i},\ s_i=(Latitude, Longitude)_i, i=36,\cdots,51\},$$ which would correspond to the measurements ozone concentration at the date of the April 13th, 2021 on these 16 others sites assumed non-visited at this same date. 
 
 \begin{table}[h!]
 \centering
        \begin{tabular}{ccc}
        \hline
        &SSFLRD &SFLRD \\\hline
        Prediction error (PE) &$0.0320$&$0.0334$\\\hline
        \end{tabular}
        \caption{Prediction error computed from both models with $h=0.32$, $\psi=0.01$, $w=0.28$.} \label{tab:2}
     \end{table}
    
    \begin{figure}[h!]
    \centering
    \includegraphics[width=0.99\textwidth]{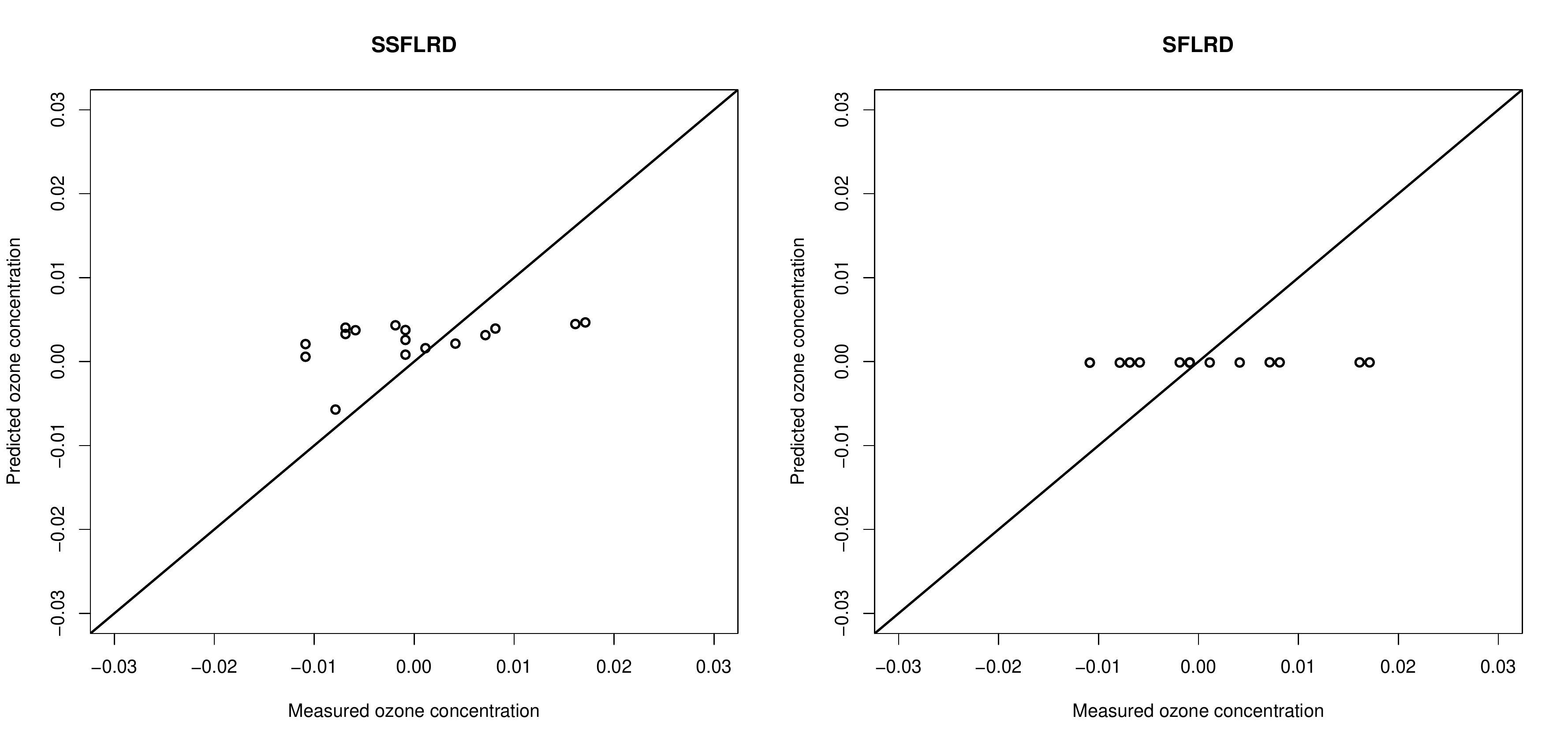}
    \caption{Centered predicted values of ozone concentration, from SSFLRD (left) and SFLRD(right), versus the centered measured values.  }\label{fig1}
    \end{figure} 

\noindent Both last graphics of Figure \ref{fig1} present a very minor different and it is confirmed by computation (see Table \ref{tab:2}) of prediction error (PE) given by
 $$PE=\sqrt{\sum_{i=36}^{51}\left(Y_{s_i}-\widehat{Y}_{s_i}\right)^2}.$$
 We see a very small advantage for prediction obtained from model SSFLRD studied in this paper.

\section{Discussion}\label{Discussion}

\noindent In this paper, we propose to study asymptotic properties of a prediction at non-visited sites computed from an estimator of nonparametric regression function in a spatial semi-functional linear regression model with derivatives. The originality of the proposed method is to consider spatially dependent data in this new model. We established the convergence rates of the estimation and prediction errors when considered processes are  $\alpha$-mixing dependent. The main contributions of this work are on convergence rates of the empirical covariance operators  and on the estimator $\widehat{r}$ of $r$ constructed from $\alpha$-mixing dependent data satisfying the general model defined in (\ref{1.2}), allowing to establish prediction at non-visited site. Its convergence rate is optimal with $d\ge4$, is better than the one of \cite{carbonetal07} with $d=3$, and is quite close from the one of \cite{carbonetal07} with $d=2$. Besides, the simulation study revealed that the presence of functional data in the SSFLRD model does not modify the convergence rate of the estimated nonparametric regression function. Notice, however, that problems are different: here it is question of prediction at non-visited sites evaluated from estimation, while \cite{carbonetal07} is only interested in that of estimation. Application to ozone pollution revealed that the proposed prediction fits well to the spatial semi-functional linear regression model with derivatives. Moreover, the SSFLRD model produces equivalent predictions with the SFLRD model. However, the presented methodology in this paper has more advantages than that of SFLRD.




\section{Proofs of asymptotic results}\label{s6} 
  
  \subsection{Technical Lemmas}
  
\noindent The proof of the following Lemma \ref{lemme1} is similar to the one of Lemma 1 in \cite{deo}.

\begin{lemma}\label{lemme1}
Assume that (\ref{ar2}) holds. Let ${\mathcal L}_r({\mathcal A})$ be the class of ${\mathcal A}$-measurable random function $X$ satisfying $\|X\|_p=\left(\E(|X|^p)\right)^{1/p}$. Let $p,\ s,\ h$ be positive constants such that $p^{-1}+s^{-1}+h^{-1}=1$, $X\in {\mathcal L}_p({\mathcal B}(S))$ and $Y\in {\mathcal L}_s({\mathcal B}(S'))$. Then
$$|\E(XY)-\E(X)\E(Y)|\leq K\left\|X\right\|_p\left\|Y\right\|_s\left\{\alpha_{1,\infty}(\rho(S,S'))\right\}^{1/h}.$$
\end{lemma}  
  
\noindent The following Lemma \ref{lemme2} adapts Proposition 8 in \cite{mas_pumo09}.

\begin{lemma}\label{lemme2}
We have: $$\left\|\left(S_{\mathbf{n},\phi}+\psi_{\mathbf{n}}I\right)^{-1}\right\|_{\infty}\le \dfrac{1}{\psi_{\mathbf{n}}},$$ where $\left\|R\right\|_{\infty}=\sup_{x\in H}\dfrac{\left\|Rx\right\|_H}{\left\|x\right\|_H}.$
\end{lemma}  

\noindent From Lemmas 9, 10 of \cite{mas_pumo09} together with Corollary 3.1 of \cite{bouka3}, we obtain the following Lemma \ref{lemme3}.

\begin{lemma}\label{lemme3}
We have:
  \begin{eqnarray*}
  \|u_{\mathbf{n},\phi}-u_{\phi}\|_{L^2({\mathcal HS})}=O\left(\frac{\log\widehat{\mathbf{n}}}{w_{\mathbf{n}}\widehat{\mathbf{n}}^{1/2}}\right) \,\,\textrm{ and }\,\,
  \|S_{\mathbf{n},\phi}-S_{\phi}\|_{L^2({\mathcal HS})}=O\left(\frac{\log\widehat{\mathbf{n}}}{w_{\mathbf{n}}\widehat{\mathbf{n}}^{1/2}}\right)\label{e2},
  \end{eqnarray*} 
  where ${\mathcal HS}$ stands for the space of Hilbert-Schmidt operators endowed with the inner product $\left\langle R, T\right\rangle_{\mathcal HS}=\sum^{+\infty}_{i=1}\left\langle R(w_i),T(w_i)\right\rangle_H$ where $w_i$ is a basis of $H$, and $\left\|R\right\|_{L^2({\mathcal HS})}=\left\{\E\left(\left\|R\right\|_{\mathcal HS}^2\right)\right\}^{1/2}.$
\end{lemma}

\begin{lemma}\label{l1}
We have:
 $$\left\|\widehat{\phi}_{\mathbf{n}}-\phi \right\|_H=O_p(1)\ \ \text{and}\ \ \left\|\widehat{\gamma}_{\mathbf{n}}-\gamma \right\|_G=O_p(1);$$
  \end{lemma}

 \noindent\textit{Proof}. 
   Putting $u_{\phi}=\Delta-\Gamma^{\prime}\Gamma^{\prime\prime-1}\Delta^{\prime}$, we have $$\widehat{\phi}_{\mathbf{n}}-\phi=(S_{\mathbf{n},\phi}+\psi_{\mathbf{n}}I)^{-1}(u_{\mathbf{n},\phi}-u_{\phi})+(S_{\mathbf{n},\phi}+\psi_{\mathbf{n}}I)^{-1}(S_{\phi}-S_{\mathbf{n},\phi}-\psi_{\mathbf{n}}I)\phi.$$ 

 \noindent From Lemma \ref{lemme2}, we obtain
  \begin{eqnarray*}
  \left\|\widehat{\phi}_{\mathbf{n}}-\phi \right\|_H^2&\leq& 3\left(\left\|u_{\mathbf{n},\phi}-u_{\phi} \right\|_H^2\left\|(S_{\mathbf{n},\phi}+\psi_{\mathbf{n}}I)^{-1} \right\|_\infty^2\right)\\
  &&+3\left\|\phi \right\|_H^2\left(\left\|S_{\phi}-S_{\mathbf{n},\phi} \right\|_\infty^2\left\|(S_{\mathbf{n},\phi}+\psi_{\mathbf{n}}I)^{-1} \right\|_\infty^2\right)\\
  &&+3\psi_{\mathbf{n}}^2\left\|\phi \right\|_H^2\left(\left\|(S_{\mathbf{n},\phi}+\psi_{\mathbf{n}}I)^{-1} \right\|_\infty^2\right)\\
  &\le&\frac{3}{\psi_{\mathbf{n}}^2}\left(\left\|u_{\mathbf{n},\phi}-u_{\phi} \right\|_H^2+\left\|S_{\phi}-S_{\mathbf{n},\phi} \right\|_\infty^2\left\|\phi\right\|_H^2\right)+3\left\|\phi\right\|_H^2.
  \end{eqnarray*}
  
\noindent For $\tau>0$, from Markov inequality, we have
\begin{eqnarray*}
&&\p\left(\left\|\widehat{\phi}_{\mathbf{n}}-\phi \right\|_H>\tau\right)\\
&=&\p\left(\left\|\widehat{\phi}_{\mathbf{n}}-\phi \right\|_H^2>\tau^2\right)\\
&\le&\p\left(\frac{3}{\psi_{\mathbf{n}}^2}\left(\left\|u_{\mathbf{n},\phi}-u_{\phi} \right\|_H^2+\left\|S_{\phi}-S_{\mathbf{n},\phi} \right\|_\infty^2\left\|\phi\right\|_H^2\right)+3\left\|\phi\right\|_H^2>\tau^2\right)\\
&\le&\p\left(\frac{3}{\psi_{\mathbf{n}}^2}\left\|u_{\mathbf{n},\phi}-u_{\phi} \right\|_H^2>\dfrac{\tau^2}{3}\right)+\p\left(\frac{3\left\|\phi\right\|_H^2}{\psi_{\mathbf{n}}^2}\left\|S_{\phi}-S_{\mathbf{n},\phi} \right\|_\infty^2>\dfrac{\tau^2}{3}\right)\\
&&+\,\p\left(3\left\|\phi\right\|_H^2>\dfrac{\tau^2}{3}\right)\\
&\le&\dfrac{9\E\left(\left\|u_{\mathbf{n},\phi}-u_{\phi} \right\|_H^2\right)}{\tau^2\psi_{\mathbf{n}}^2}+\dfrac{9\left\|\phi\right\|_H^2\E\left(\left\|S_{\phi}-S_{\mathbf{n},\phi} \right\|_\infty^2\right)}{\tau^2\psi_{\mathbf{n}}^2}+\mathbb{I}_{]0,3\|\phi\|_H[}(\tau),
\end{eqnarray*}
where $\mathbb{I}_{]0,3\|\phi\|_H[}(.)$ is the indicator function. Finally, from Lemma \ref{lemme3}, we conclude that $$\forall\delta>0,\, \exists\tau\ge3\|\phi\|_H,\,N_{\delta}\in\N\, \, \text{such that}\, \, \forall\mathbf{n}>N_{\delta}\mathbf{1},\, \p\left(\left\|\widehat{\phi}_{\mathbf{n}}-\phi \right\|_H>\tau\right)\le\delta,$$
where $\mathbf{n}>N_{\delta}\mathbf{1}$ means that $\min_{1\le k\le d}\{n_k\}>N_{\delta}$. Similarly, we obtain that $\left\|\widehat{\gamma}_{\mathbf{n}}-\gamma \right\|_G=O_p(1)$.

  \subsection{Proof of Theorem \ref{th1}}
  
 \noindent Proof of  (\ref{3.2}):  For seek of simplicity, we only give the proof of the empirical operator $\Gamma_{\mathbf{n}}$, since for other estimators similar arguments can be applied. Recall that $\Lambda_{\mathbf{n}}=\Gamma_{\mathbf{n}}$ with $\Gamma_{\mathbf{n}}=\dfrac{1}{\widehat{\mathbf{n}}}\sum_{\mathbf{i}\in{\cal I}_{\mathbf{n}}}X_{\mathbf{i}}\otimes_{H}X_{\mathbf{i}}$ and $\Lambda=\Gamma=\E(X\otimes_{H}X)$. For all $\tau>0$ and From Markov inequality,  we have:
  \begin{eqnarray}
  \p\left(\left\|\Gamma_{\mathbf{n}}-\Gamma\right\|_{\infty}>\tau\right)&=&\p\left(\sup_{k}\left\|\Gamma_{\mathbf{n}}v_k-\Gamma v_k\right\|_{H}>\tau\right)\nonumber\\
  &\leq&\sum^{+\infty}_{k=1}\p\left(\left\|\Gamma_{\mathbf{n}}v_k-\Gamma v_k\right\|_{H}>\tau\right)\nonumber\\
  &\leq&\dfrac{1}{\tau^2}\sum^{\infty}_{k=1}\E\left(\left\|\Gamma_{\mathbf{n}}v_k-\Gamma v_k\right\|_{H}^2\right)\label{6.1}.
  \end{eqnarray}
However, Setting 
  \begin{gather*}
  L_{\mathbf{ij}}^k=\left\langle\left\langle X_{\mathbf{i}}, v_k\right\rangle_{H}X_{\mathbf{i}}-\E(\left\langle X, v_k\right\rangle_{H}X),\left\langle X_{\mathbf{j}}, v_k\right\rangle_{H}X_{\mathbf{j}}-\E(\left\langle X, v_k\right\rangle_{H}X) \right\rangle_{H},
  \end{gather*}
   we have
  \begin{eqnarray*}
  \E\left[\|\Gamma_{\mathbf{n}}v_k-\Gamma v_k\|^{2}_{H}\right]&=&\E\left[\left\|\frac{1}{\widehat{\mathbf{n}}}\sum_{\mathbf{i}\in{\cal I}_{\mathbf{n}}}\left(\left\langle X_{\mathbf{i}}, v_k\right\rangle_{H}X_{\mathbf{i}}-\E(\left\langle X, v_k\right\rangle_{H}X)\right)\right\|^{2}_{H}\right]\\
  &=&\dfrac{1}{\widehat{\mathbf{n}}^{2}}\sum_{\mathbf{i}\in{\cal I}_{\mathbf{n}}}\E\left[\left\|\left\langle X_{\mathbf{i}}, v_k\right\rangle_{H}X_{\mathbf{i}}-\E(\left\langle X, v_k\right\rangle_{H}X)\right\|^{2}_{H}\right]\\
  &&+\dfrac{1}{\widehat{\mathbf{n}}^{2}}\sum_{\mathbf{i}\neq\mathbf{j}}\E \left(L_{\mathbf{ij}}^k\right)\\
  &:=&A_k+B_k.
  \end{eqnarray*}
  On the one hand, since $X_{\mathbf{i}}$ are strictly stationary with the same law as $X$, from Assumption \ref{as5}, we have:
  \begin{eqnarray}
  A_k&=&\dfrac{1}{\widehat{\mathbf{n}}^{2}}\sum_{\mathbf{i}\in{\cal I}_{\mathbf{n}}}\E\left[\left\|\left\langle X_{\mathbf{i}}, v_k\right\rangle_{H}X_{\mathbf{i}}-\E(\left\langle X, v_k\right\rangle_{H}X)\right\|^{2}_{H}\right]\nonumber\\
  &\leq&\dfrac{2}{\widehat{\mathbf{n}}^{2}}\sum_{\mathbf{i}\in{\cal I}_{\mathbf{n}}}\left(\E\left(\left\langle X_{\mathbf{i}}, v_k\right\rangle^{2}_{H}\|X_{\mathbf{i}}\|^{2}_{H}\right)+\E\left(\left\langle X, v_k\right\rangle^{2}_{H}\|X\|^{2}_{H}\right)\right)\nonumber\\
  &\leq&\dfrac{4C^2\lambda_k}{\widehat{\mathbf{n}}}.\label{6.2}
  \end{eqnarray}
  On the other hand, we have: 
  \begin{eqnarray*}
  B_k&=&\dfrac{1}{\widehat{\mathbf{n}}^{2}}\sum_{0<\|\mathbf{i}-\mathbf{j}\|\leq C_{\mathbf{n}}}\E (L_{\mathbf{ij}}^k)+\dfrac{1}{\widehat{\mathbf{n}}^{2}}\sum_{\|\mathbf{i}-\mathbf{j}\|>C_{\mathbf{n}}}\E (L_{\mathbf{ij}}^k):=B_{k1}+B_{k2}
  \end{eqnarray*}
  where $0<C_{\mathbf{n}}<\widehat{\mathbf{n}}$ and $C_{\mathbf{n}}\rightarrow+\infty$ as $\mathbf{n}\rightarrow+\infty$. However, from Cauchy-Schwartz inequality and by the same arguments than those of the relation (\ref{6.2}), we have 
  \begin{eqnarray*}
  \E(|L_{\mathbf{ij}}^k|)\leq\E(\sqrt{L_{\mathbf{ii}}^k}\sqrt{L_{\mathbf{jj}}^k})\leq (\E (L_{\mathbf{ii}}^k))^{1/2}(\E (L_{\mathbf{jj}}^k))^{1/2}\leq 2C^2\lambda_k 
  \end{eqnarray*}
  because $\E (L_{\mathbf{ii}}^k)=\E\left(\left\|\left\langle X_{\mathbf{i}}, v_k\right\rangle_{H}X_{\mathbf{i}}-\E\left(\left\langle X, v_k\right\rangle_{H}X\right)\right\|^{2}_{H}\right)$ are terms of $A_k$. Then
  \begin{eqnarray*}
  |B_{k1}|&\leq&\dfrac{2C^2\lambda_k}{\widehat{\mathbf{n}}^2}\sum^{C_{\mathbf{n}}}_{\ell=1}\sum_{\stackrel{\mathbf{i,j}\in{\mathcal I}_{\mathbf{n}}}{\ell\leq\|\mathbf{i-j}\|=t<\ell+1}}1\\
  &\leq&\dfrac{2C^2\lambda_k}{\widehat{\mathbf{n}}^2}\sum^{C_{\mathbf{n}}}_{\ell=1}\sum_{\mathbf{j}\in{\mathcal I}_{\mathbf{n}}}\sum_{\stackrel{\mathbf{i}\in{\mathcal I}_{\mathbf{n}}}{\ell\leq\|\mathbf{i-j}\|=t<\ell+1}}1\\
  &\leq&\dfrac{2C^2\lambda_k}{\widehat{\mathbf{n}}}\sum^{C_{\mathbf{n}}}_{\ell=1}\sum_{\stackrel{\mathbf{i}\in{\mathcal I}_{\mathbf{n}}}{\ell\leq\|\mathbf{i}\|=t<\ell+1}}1\\
  &\leq&\dfrac{2C^2\lambda_k}{\widehat{\mathbf{n}}}\sum^{C_{\mathbf{n}}}_{t=1}t^{d-1}\leq\dfrac{2C^2\lambda_kC_{\mathbf{n}}^d}{\widehat{\mathbf{n}}}.
  \end{eqnarray*}
  Taking $C_{\mathbf{n}}=\lfloor(\log\widehat{\mathbf{n}})^{1/d}\rfloor$ (where $\lfloor x\rfloor$ stands for the integer part of $x$), we obtain
  \begin{eqnarray}\label{6.3}
  |B_{k1}|\leq\dfrac{2C^2\lambda_k\log\widehat{\mathbf{n}}}{\widehat{\mathbf{n}}}
  \end{eqnarray}
  Since $X_{\mathbf{i}}$ are strictly stationary with the same law as $X$ and $\|X\|_{H}<C$ a.s., it follows that $$\E\left[\left(L^k_{\mathbf{ii}}\right)^4\right]=\E\left(\|\left\langle X_{\mathbf{i}}, v_k\right\rangle_{H}X_{\mathbf{i}}-\E\left(\left\langle X, v_k\right\rangle_{H}X\right)\|_{H}^4\right)\leq8\E\left(\left\langle X,v_k\right\rangle_{H}^4\|X\|_{H}^4\right)<8\lambda_{k}C^{6}.$$ Applying Lemma \ref{lemme1}  with $p=s=4$ and $h=2$, we obtain
  \begin{eqnarray*}
  |\E (L_{\mathbf{ij}}^k)|\leq K\|L_{\mathbf{ii}}^k\|^2_{4}\{\alpha_{1,\infty}(|\mathbf{i-j}|)\}^{1/2}\leq 2KC^3\sqrt{2\lambda_k}\{\alpha_{1,\infty}(|\mathbf{i-j}|)\}^{1/2}.
  \end{eqnarray*}
  Since $\alpha_{1,\infty}(t)=O(t^{-\theta})$ with $\theta>2d$, then
  \begin{eqnarray}
  |B_{k2}|&\leq&\dfrac{2KC^3\sqrt{2\lambda_k}}{\widehat{\mathbf{n}}^{2}}\sum^{+\infty}_{\ell=C_{\mathbf{n}}+1}\sum_{\stackrel{\mathbf{i,j}\in{\mathcal I}_{\mathbf{n}}}{\ell\leq\|\mathbf{i-j}\|=t<\ell+1}}\{\alpha_{1,\infty}(t)\}^{1/2}\nonumber\\
  &\leq&\dfrac{2KC^3\sqrt{2\lambda_k}}{\widehat{\mathbf{n}}}\sum^{+\infty}_{\ell=C_{\mathbf{n}}+1}\sum_{\stackrel{\mathbf{i}\in{\mathcal I}_{\mathbf{n}}}{\ell\leq\|\mathbf{i}\|=t<\ell+1}}\{\alpha_{1,\infty}(t)\}^{1/2}\nonumber\\
  &\leq&\dfrac{2KC^3\sqrt{2\lambda_k}}{\widehat{\mathbf{n}}}\sum^{+\infty}_{t=1}t^{d-1}\{\alpha_{1,\infty}(t)\}^{1/2}\nonumber\\
  &\leq&\dfrac{2KC^3\sqrt{2\lambda_k}}{\widehat{\mathbf{n}}}\sum^{\infty}_{t=1}t^{d-1-\theta/2}.\label{6.4}
  \end{eqnarray}
  From inequalities (\ref{6.1}), (\ref{6.2}), (\ref{6.3}) and (\ref{6.4}) with $\lambda_k=O(u^k)$, $0<u<1$, we conclude that
  \begin{eqnarray*}
  \p\left(\left\|\Gamma_{\mathbf{n}}-\Gamma\right\|_{\infty}>\tau\right)=O\left(\dfrac{\log\widehat{\mathbf{n}}}{\widehat{\mathbf{n}}}\right).
  \end{eqnarray*}

 \noindent Proof of (\ref{3.3}): Take $\Lambda_{\mathbf{n}}=\Gamma_{\mathbf{n}}$ with $\Gamma_{\mathbf{n}}=\dfrac{1}{\widehat{\mathbf{n}}}\sum_{\mathbf{i}\in{\cal I}_{\mathbf{n}}}X_{\mathbf{i}}\otimes_{H}X_{\mathbf{i}}$ and $\Lambda=\Gamma=\E(X\otimes_{H}X)$. By definition, we have
  $$
  \left\|\Lambda_{\mathbf{n}}-\Lambda\right\|_{L^{2}({\cal HS})}=\left\{\E\left[\left\|\Lambda_{\mathbf{n}}-\Lambda\right\|^{2}_{{\cal HS}}\right]\right\}^{1/2}.
  $$
  Since $(v_{j})_{j\geq1}$ is an orthonormal basis of $H$, we have
  \begin{eqnarray*}
  &&\E\left[\left\|\Lambda_{\mathbf{n}}-\Lambda\right\|^{2}_{{\cal HS}}\right]=\sum^{+\infty}_{i=1}\E\left[\left\|\Lambda_{\mathbf{n}}(v_{i})-\Lambda(v_{i})\right\|^{2}_{H}\right]\\
  &&\ \ \ =\sum^{Q}_{i=1}\E\left[\left\|\Lambda_{\mathbf{n}}(v_{i})-\Lambda(v_{i})\right\|^{2}_{H}\right]+\sum_{i>Q}\E\left[\left\|\Lambda_{\mathbf{n}}(v_{i})-\Lambda(v_{i})\right\|^{2}_{H}\right]:=A+B.
  \end{eqnarray*}
  Applying Theorem \ref{th1} and taking $Q=\lfloor K\log\widehat{\mathbf{n}}\rfloor$, we have
  $$
  A\leq\E\left[\left\|\Lambda_{\mathbf{n}}-\Lambda\right\|^{2}_{\infty}\right]\sum^{Q}_{i=1}\left\|v_{i}\right\|^{2}_{H}\leq C_1\dfrac{(\log\widehat{\mathbf{n}})^{2}}{\widehat{\mathbf{n}}}
  $$
  where $C_1$ is some positive constant. On the other hand, we have
  \begin{eqnarray*}
  B&=&\sum_{i>Q}\E\left[\sum^{+\infty}_{j=1}\left\langle \Lambda_{\mathbf{n}}(v_{i})-\Lambda(v_{i}), v_{j}\right\rangle^{2}_{H}\right]\\
  &=&\sum_{i>Q}\sum^{+\infty}_{j=1}\E\left\{\left[\frac{1}{\widehat{\mathbf{n}}}\sum_{\mathbf{k}\in{\cal I}_{\mathbf{n}}}\left\langle \left\langle X_{\mathbf{k}}, v_{i}\right\rangle_{H}X_{\mathbf{k}}-\E(\left\langle X, v_{i}\right\rangle_{H}X), v_{j}\right\rangle_{H}\right]^{2}\right\}\\
  &\leq&\dfrac{1}{\widehat{\mathbf{n}}}\sum_{i>Q}\sum^{+\infty}_{j=1}\sum_{\mathbf{k}\in{\cal I}_{\mathbf{n}}}\E\left[\left\langle \left\langle X_{\mathbf{k}}, v_{i}\right\rangle_{H}X_{\mathbf{k}}-\E(\left\langle X, v_{i}\right\rangle_{H}X), v_{j}\right\rangle^{2}_{H}\right]\\
  &\leq&\frac{2}{\widehat{\mathbf{n}}}\sum_{i>Q}\sum^{+\infty}_{j=1}\sum_{\mathbf{k}\in{\cal I}_{\mathbf{n}}}\left[\E\left(\left\langle X_{\mathbf{k}}, v_{i}\right\rangle^{2}_{H}\left\langle X_{\mathbf{k}}, v_{j}\right\rangle^{2}_{H}\right)+\E\left(\left\langle X, v_{i}\right\rangle^{2}_{H}\left\langle X, v_{j}\right\rangle^{2}_{H}\right)\right]\\
  &\leq&4\sum_{i>Q}\left[\E\left(\left\langle X, v_{i}\right\rangle^{4}_{H}\right)\right]^{1/2}\sum^{+\infty}_{j=1}\left[\E\left(\left\langle X, v_{j}\right\rangle^{4}_{H}\right)\right]^{1/2}.
  \end{eqnarray*}
  Since $\left\langle X, v_{j}\right\rangle^{4}_{H}\leq\|X\|^{2}_{H}\|v_{j}\|^{2}_{H}\left\langle X, v_{j}\right\rangle^{2}_{H}<C^{2}\left\langle X, v_{j}\right\rangle^{2}_{H}$ a.s. and $\E\left(\left\langle X, v_{j}\right\rangle^{2}_{H}\right)=\lambda_{j}$ with $\lambda_{j}=O(u^{j})$, $0<u<1$, $j\geq1$, and $Q=\lfloor K\log\widehat{\mathbf{n}}\rfloor$ with $K=\dfrac{3}{\log\left(\dfrac{1}{u}\right)}$, then
  $$
  B\leq C_{2}\exp(-K(\log\widehat{\mathbf{n}})(\log(1/u))/2)=\dfrac{C_{2}}{\widehat{\mathbf{n}}^{3/2}}
  $$
  where $C_{2}$ is a positive constant. This finishes the proof of  $(\ref{3.3})$.

 \subsection{Proof of Theorem \ref{regl2}}

\noindent  Put  $K_{\mathbf{i}}:=\dfrac{1}{h^{d}}K\left(\dfrac{\dfrac{\mathbf{i}}{\mathbf{n+1}}-\mathbf{s}_{0}}{h}\right)$, $\Gamma(\mathbf{i}):=\left(\dfrac{\mathbf{i}}{\mathbf{n+1}}-\mathbf{s}_{0}\right)^{T}r''(\mathbf{s}_{0})\left(\dfrac{\mathbf{i}}{\mathbf{n+1}}-\mathbf{s}_{0}\right)$, where $r''(\mathbf{s}_{0})$ stands for the matrix of second order partial derivatives of $r$ at $\mathbf{s}_{0}$,
 \begin{eqnarray*}
 A_{\mathbf{n}}&=&\frac{1}{\widehat{\mathbf{n}}}{\mathcal X}^{T}W_{0}{\mathcal X}\\
 &=&
 \begin{pmatrix}
 	\dfrac{1}{\widehat{\mathbf{n}}}\sum_{\mathbf{i}\in{\mathcal I}_{\mathbf{n}}}K_{\mathbf{i}}&\dfrac{1}{\widehat{\mathbf{n}}}\sum_{\mathbf{i}\in{\mathcal I}_{\mathbf{n}}}\left(\frac{\frac{\mathbf{i}}{\mathbf{n+1}}-\mathbf{s}_{0}}{h}\right)^{T}K_{\mathbf{i}}\\
 	\dfrac{1}{\widehat{\mathbf{n}}}\sum_{\mathbf{i}\in{\mathcal I}_{\mathbf{n}}}\left(\frac{\frac{\mathbf{i}}{\mathbf{n+1}}-\mathbf{s}_{0}}{h}\right)K_{\mathbf{i}}&\ \ \ \ \ \ \dfrac{1}{\widehat{\mathbf{n}}}\sum_{\mathbf{i}\in{\mathcal I}_{\mathbf{n}}}\left(\frac{\frac{\mathbf{i}}{\mathbf{n+1}}-\mathbf{s}_{0}}{h}\right)\left(\frac{\frac{\mathbf{i}}{\mathbf{n+1}}-\mathbf{s}_{0}}{h}\right)^{T}K_{\mathbf{i}}\\
 \end{pmatrix}
 \end{eqnarray*}
 and
 \[B_{\mathbf{n}}:=\frac{1}{\widehat{\mathbf{n}}}{\mathcal X}^{T}W_{0}{\mathcal Y}=
 \begin{pmatrix}
 	\dfrac{1}{\widehat{\mathbf{n}}}\sum_{\mathbf{i}\in{\mathcal I}_{\mathbf{n}}}K_{\mathbf{i}}T_{\mathbf{i}}\\
 	\dfrac{1}{\widehat{\mathbf{n}}}\sum_{\mathbf{i}\in{\mathcal I}_{\mathbf{n}}}\left(\dfrac{\dfrac{\mathbf{i}}{\mathbf{n+1}}-\mathbf{s}_{0}}{h}\right)K_{\mathbf{i}}T_{\mathbf{i}}\\
 \end{pmatrix},\ W_{\mathbf{n}}:=B_{\mathbf{n}}-A_{\mathbf{n}}\begin{pmatrix}
 	\beta_{0}\\
 	h\beta_{1}\\
 \end{pmatrix}.
 \]
 We have
 \[\lim_{\mathbf{n}\to +\infty}A_{\mathbf{n}}=
 \begin{pmatrix}
  1&\mathbf{0}^{T}\\
 	\mathbf{0}&\ \ \ \int \mathbf{u}\mathbf{u}^{T}K(\mathbf{u})d\mathbf{u}\\	
 \end{pmatrix}=\begin{pmatrix}
  1&\mathbf{0}^{T}\\
 	\mathbf{0}&\ \ \ \nu_2I_d\\	
 \end{pmatrix}=A,  
 \]
 and $ \textrm{det}(A)=\nu_2 \ne 0$;  then $A$ is invertible and
 \[A^{-1}= 
 \begin{pmatrix}
  1&\mathbf{0}^{T}\\
 	\mathbf{0}&\ \ \ \nu_2^{-1}I_d\\	
 \end{pmatrix}.
 \]
 Therefore, putting
 \[
 A^{-1}_{\mathbf{n}} =\begin{pmatrix}
 	u^{\mathbf{n}}_{11}&u^{\mathbf{n}}_{12}\\
 	u^{\mathbf{n}}_{21}&u^{\mathbf{n}}_{22}\\
 \end{pmatrix},
 \]
 we have $\lim_{\mathbf{n}\rightarrow +\infty}\left(u^{\mathbf{n}}_{11}\right)=1$ and $\lim_{\mathbf{n}\rightarrow +\infty}\left(u^{\mathbf{n}}_{12}\right)=\mathbf{0}^{T}$,  and the sequences $\left(u^{\mathbf{n}}_{11}\right)_{\mathbf{n}}$ and  $\left(u^{\mathbf{n}}_{12}\right)_{\mathbf{n}}$ are bounded.
 
 \noindent $(i)$ For all $\mathbf{s}_{0}\in[0,1]^d$, we have
 $
 \widehat{r}(\mathbf{s}_{0})-r(\mathbf{s}_{0})=(1,\mathbf{0}^{T})A^{-1}_{\mathbf{n}}W_{\mathbf{n}}$ and
 \[\E\left(\widehat{r}(\mathbf{s}_{0})-r(\mathbf{s}_{0})\right)=(1,\mathbf{0}^{T})A^{-1}_{\mathbf{n}}
 \begin{pmatrix}
 \dfrac{1}{\widehat{\mathbf{n}}}\sum_{\mathbf{i}\in{\mathcal I}_{\mathbf{n}}}K_{\mathbf{i}}(\Gamma(\mathbf{i})+\E(\epsilon^*_{\mathbf{i}}))\\
 	\dfrac{1}{\widehat{\mathbf{n}}}\sum_{\mathbf{i}\in{\mathcal I}_{\mathbf{n}}}\left(\dfrac{\dfrac{\mathbf{i}}{\mathbf{n+1}}-\mathbf{s}_{0}}{h}\right)K_{\mathbf{i}}(\Gamma(\mathbf{i})+\E(\epsilon^*_{\mathbf{i}}))\\	
 \end{pmatrix},
 \]
 where 
 $
 \E(\epsilon^*_{\mathbf{i}})=\E\left(\left\langle \phi-\widehat{\phi}_{\mathbf{n}}, X_{\mathbf{i}}\right\rangle_{H}\right)+\E\left(\left\langle \gamma-\widehat{\gamma}_{\mathbf{n}}, X^{\prime}_{\mathbf{i}}\right\rangle_{G}\right)$. Since $$\E\left(\left\langle(S_{\mathbf{n},\phi}+\psi_{\mathbf{n}}I)^{-1}\phi, X\right\rangle_H\right)\to\E\left(\left\langle S_{\phi}^{-1}\phi, X\right\rangle_H\right)<+\infty,$$
 from Assumption \ref{as5} and from Lemmas \ref{lemme2} and \ref{lemme3}, we have
 \begin{eqnarray*}
  \E\left(\left\langle\widehat{\phi}_{\mathbf{n}}-\phi,X\right\rangle_H\right)&=&\E\left(\left\langle(S_{\mathbf{n},\phi}+\psi_{\mathbf{n}}I)^{-1}(u_{\mathbf{n},\phi}-u_{\phi}),X\right\rangle_H\right)\\
  &&-\E\left(\left\langle(S_{\mathbf{n},\phi}+\psi_{\mathbf{n}}I)^{-1}\phi,\psi_{\mathbf{n}} X\right\rangle_H\right)\\
 &&  +\E\left(\left\langle(S_{\mathbf{n},\phi}+\psi_{\mathbf{n}}I)^{-1}(S_{\phi}-S_{\mathbf{n},\phi})\phi, X\right\rangle_H\right)\\
 &=& O\left(\dfrac{\log \widehat{\mathbf{n}}}{w_{\mathbf{n}}\psi_{\mathbf{n}}\widehat{\mathbf{n}}^{1/2}}\right)+O\left(\psi_{\mathbf{n}}\right).
 \end{eqnarray*}
 Therefore,
  $$\left(\dfrac{1}{\widehat{\mathbf{n}}}\sum_{\mathbf{i}\in{\mathcal I}_{\mathbf{n}}}K_{\mathbf{i}}\E\left[\left(\left\langle \phi-\widehat{\phi}_{\mathbf{n}}, X_{\mathbf{i}}\right\rangle_{H}\right)\right]\right)^2=O\left(\dfrac{(\log \widehat{\mathbf{n}})^2}{w_{\mathbf{n}}^2\psi_{\mathbf{n}}^2\widehat{\mathbf{n}}}\right)+O\left(\psi_{\mathbf{n}}^2\right).$$
 Similarly, we have $$\left(\dfrac{1}{\widehat{\mathbf{n}}}\sum_{\mathbf{i}\in{\mathcal I}_{\mathbf{n}}}K_{\mathbf{i}}\E\left[\left\langle \gamma-\widehat{\gamma}_{\mathbf{n}}, X^{\prime}_{\mathbf{i}}\right\rangle_{G}\right]\right)^2=O\left(\dfrac{(\log \widehat{\mathbf{n}})^2}{w_{\mathbf{n}}^2\psi_{\mathbf{n}}^2\widehat{\mathbf{n}}}\right)+O\left(\psi_{\mathbf{n}}^2\right)$$  and $$	\left(\dfrac{1}{\widehat{\mathbf{n}}}\sum_{\mathbf{i}\in{\mathcal I}_{\mathbf{n}}}\left(\frac{\dfrac{\mathbf{i}}{\mathbf{n+1}}-\mathbf{s}_{0}}{h}\right)K_{\mathbf{i}}\E(\epsilon^*_{\mathbf{i}})\right)^2=O\left(\dfrac{(\log \widehat{\mathbf{n}})^2}{w_{\mathbf{n}}^2\psi_{\mathbf{n}}^2\widehat{\mathbf{n}}}\right)+O\left(\psi_{\mathbf{n}}^2\right).$$
 Besides, we have $$\dfrac{1}{\widehat{\mathbf{n}}}\sum_{\mathbf{i}\in{\mathcal I}_{\mathbf{n}}}K_{\mathbf{i}}\Gamma(\mathbf{i})=O(h^2)\, \, \,  \text{and}\, \, \,  \dfrac{1}{\widehat{\mathbf{n}}}\sum_{\mathbf{i}\in{\mathcal I}_{\mathbf{n}}}\left(\dfrac{\dfrac{\mathbf{i}}{\mathbf{n+1}}-\mathbf{s}_{0}}{h}\right)K_{\mathbf{i}}\Gamma(\mathbf{i})=O(h^2).$$ Thus $$\sup_{\mathbf{s}_0\in[0,1]^d}\bigg\{\E\bigg(\widehat{r}(\mathbf{s}_{0})-r(\mathbf{s}_{0})\bigg)\bigg\}^2=O\left(h^4\right)+O\left(\dfrac{(\log \widehat{\mathbf{n}})^2}{w_{\mathbf{n}}^2\psi_{\mathbf{n}}^2\widehat{\mathbf{n}}}\right)+O\left(\psi_{\mathbf{n}}^2\right).$$
 
 \smallskip
\noindent $(ii)$    Since $T_{\mathbf{i}}-\beta_{0}-\left(\dfrac{\dfrac{\mathbf{i}}{\mathbf{n+1}}-\mathbf{s}_{0}}{h}\right)^{T}h\beta_{1}-\Gamma(\mathbf{i})=\xi_{\mathbf{i}}$ and putting $k_{\mathbf{n}}(\mathbf{i})=u^{\mathbf{n}}_{11}+u^{\mathbf{n}}_{12}\left(\dfrac{\dfrac{\mathbf{i}}{\mathbf{n+1}}-\mathbf{s}_{0}}{h}\right)$, we obtain for all $\mathbf{s}_0\in[0,1]^d$ that
 \begin{eqnarray*}
 \widehat{r}(\mathbf{s}_{0})-\E\left(\widehat{r}(\mathbf{s}_{0})\right)&=&\dfrac{1}{\widehat{\mathbf{n}}}\sum_{\mathbf{i}\in{\mathcal I}_{\mathbf{n}}}K_{\mathbf{i}}\left(\xi_{\mathbf{i}}-\E\left(\xi_{\mathbf{i}}\right)\right)k_{\mathbf{n}}(\mathbf{i})\\
 &=&\dfrac{1}{\widehat{\mathbf{n}}h^{d}}\sum_{\mathbf{i}\in{\mathcal I}_{\mathbf{n}}}K\left(\dfrac{\dfrac{\mathbf{i}}{\mathbf{n+1}}-\mathbf{s}_{0}}{h}\right)\left(\xi_{\mathbf{i}}-\E\left(\xi_{\mathbf{i}}\right)\right)k_{\mathbf{n}}(\mathbf{i}).
 \end{eqnarray*}
 Also, since from Lemma \ref{l1}, we have $\E(\xi_{\mathbf{i}}^2)=O(1)$, it follows that
 \begin{eqnarray*}
 \E\bigg(\bigg(\hat{r}(\mathbf{s}_{0})-\E\left(\hat{r}(\mathbf{s}_{0})\right)\bigg)^{2}\bigg)\leq\dfrac{C_1}{\left(\widehat{\mathbf{n}}h^d\right)^{2}}\sum_{\mathbf{i}\in{\mathcal I}_{\mathbf{n}}}K\left(\dfrac{\dfrac{\mathbf{i}}{\mathbf{n+1}}-\mathbf{s}_{0}}{h}\right)^{2}k^{2}_{\mathbf{n}}(\mathbf{i})+F ,
 \end{eqnarray*}
  where $C_1$ is some positive constant and $$F=\dfrac{1}{\left(\widehat{\mathbf{n}}h^d\right)^{2}}\sum_{\mathbf{i}\ne\mathbf{j}}K\left(\dfrac{\dfrac{\mathbf{i}}{\mathbf{n+1}}-\mathbf{s}_{0}}{h}\right)K\left(\dfrac{\dfrac{\mathbf{j}}{\mathbf{n+1}}-\mathbf{s}_{0}}{h}\right)k_{\mathbf{n}}(\mathbf{i})k_{\mathbf{n}}(\mathbf{j})Cov\left(\xi_{\mathbf{i}},\xi_{\mathbf{j}}\right).$$ 
  Since the sequences $\left(u^{\mathbf{n}}_{11}\right)_{\mathbf{n}}$ and  $\left(u^{\mathbf{n}}_{12}\right)_{\mathbf{n}}$ are bounded, $\left|\left(\dfrac{\dfrac{\mathbf{i}}{\mathbf{n+1}}-\mathbf{s}_{0}}{h}\right)\right|\le\mathbf{1}$, and  $$Cov\left(\xi_{\mathbf{i}},\xi_{\mathbf{j}}\right)=Cov\left(\epsilon_{\mathbf{i}},\epsilon_{\mathbf{j}}\right)+Cov\left(\epsilon^*_{\mathbf{i}},\epsilon^*_{\mathbf{j}}\right),$$
  where $\epsilon^*_{\mathbf{i}}=\left\langle \phi-\widehat{\phi}_{\mathbf{n}}, X_{\mathbf{i}}\right\rangle_{H}+\left\langle \gamma-\widehat{\gamma}_{\mathbf{n}}, X^{\prime}_{\mathbf{i}}\right\rangle_{G},$ it follows that $F=F_1+F_2$, where
 \begin{eqnarray*}
 F_1&=&\dfrac{1}{\left(\widehat{\mathbf{n}}h^d\right)^{2}}\sum_{\mathbf{i}\ne\mathbf{j}}K\left(\dfrac{\dfrac{\mathbf{i}}{\mathbf{n+1}}-\mathbf{s}_{0}}{h}\right)K\left(\dfrac{\dfrac{\mathbf{j}}{\mathbf{n+1}}-\mathbf{s}_{0}}{h}\right)k_{\mathbf{n}}(\mathbf{i})k_{\mathbf{n}}(\mathbf{j})Cov\left(\epsilon_{\mathbf{i}},\epsilon_{\mathbf{j}}\right);\\
 F_2&=&\dfrac{1}{\left(\widehat{\mathbf{n}}h^d\right)^{2}}\sum_{\mathbf{i}\ne\mathbf{j}}K\left(\dfrac{\dfrac{\mathbf{i}}{\mathbf{n+1}}-\mathbf{s}_{0}}{h}\right)K\left(\dfrac{\dfrac{\mathbf{j}}{\mathbf{n+1}}-\mathbf{s}_{0}}{h}\right)k_{\mathbf{n}}(\mathbf{i})k_{\mathbf{n}}(\mathbf{j})Cov\left(\epsilon_{\mathbf{i}}^*,\epsilon_{\mathbf{j}}^*\right). 
 \end{eqnarray*}
 Under Assumption \ref{rega2}, we have
 \begin{eqnarray*}
 &&\widehat{\mathbf{n}}h^dF_1\\
 &\le&\dfrac{c}{\widehat{\mathbf{n}}h^d}\sum_{\mathbf{i}\ne\mathbf{j}}\left|K\left(\dfrac{\dfrac{\mathbf{i}}{\mathbf{n+1}}-\mathbf{s}_{0}}{h}\right)-K\left(\dfrac{\dfrac{\mathbf{j}}{\mathbf{n+1}}-\mathbf{s}_{0}}{h}\right)\right|K\left(\dfrac{\dfrac{\mathbf{j}}{\mathbf{n+1}}-\mathbf{s}_{0}}{h}\right)\\
 &&\ \ \ \ \times\  |Cov\left(\epsilon_{\mathbf{i}},\epsilon_{\mathbf{j}}\right)|+\dfrac{c}{\widehat{\mathbf{n}}h^d}\sum_{\mathbf{i}\ne\mathbf{j}}\left[K\left(\dfrac{\dfrac{\mathbf{j}}{\mathbf{n+1}}-\mathbf{s}_{0}}{h}\right)\right]^2|Cov\left(\epsilon_{\mathbf{i}},\epsilon_{\mathbf{j}}\right)|\\
 &\le&\dfrac{c\sigma^2}{\widehat{\mathbf{n}}h^d}\sum_{\mathbf{j}\in{\mathcal I}_{\mathbf{n}}}K\left(\dfrac{\dfrac{\mathbf{j}}{\mathbf{n+1}}-\mathbf{s}_{0}}{h}\right)\sum_{\stackrel{\mathbf{i}\in{\mathcal I}_{\mathbf{n}}}{\|\mathbf{i-j}\|>0}}\dfrac{\left\|\dfrac{\mathbf{i}}{\mathbf{n+1}}-\dfrac{\mathbf{j}}{\mathbf{n+1}}\right\|}{h}\exp(-a\|\mathbf{i-j}\|)\\
 &&+\dfrac{c\sigma^2}{\widehat{\mathbf{n}}h^d}\sum_{\mathbf{j}\in{\mathcal I}_{\mathbf{n}}}\left[K\left(\dfrac{\dfrac{\mathbf{j}}{\mathbf{n+1}}-\mathbf{s}_{0}}{h}\right)\right]^2\sum_{\stackrel{\mathbf{i}\in{\mathcal I}_{\mathbf{n}}}{\|\mathbf{i-j}\|>0}}\exp(-a\|\mathbf{i-j}\|)\\
 &\leq&\dfrac{c\sigma^2}{\widehat{\mathbf{n}}h^d\min_{k=1,\cdots,d}(n_k)h}\sum_{\mathbf{j}\in{\mathcal I}_{\mathbf{n}}}K\left(\dfrac{\dfrac{\mathbf{j}}{\mathbf{n+1}}-\mathbf{s}_{0}}{h}\right)\sum^{+\infty}_{t=1}t^de^{-at}\\
 &&+\dfrac{c\sigma^2}{\widehat{\mathbf{n}}h^d}\sum_{\mathbf{j}\in{\mathcal I}_{\mathbf{n}}}\left[K\left(\dfrac{\dfrac{\mathbf{j}}{\mathbf{n+1}}-\mathbf{s}_{0}}{h}\right)\right]^2\sum^{+\infty}_{t=1}t^{d-1}e^{-at}\\
 &&\longrightarrow c_1\sigma^2\int K^2(\mathbf{u})d\mathbf{u},
 \end{eqnarray*}
 where $c$ and $c_1$ are positive constants. On the other hand, taking $Q=\lfloor (\log\widehat{\mathbf{n}})^{1/d} \rfloor$ and applying Lemma \ref{lemme1} together with Lemma \ref{l1}, we have
 \begin{eqnarray*}
 &&\widehat{\mathbf{n}}h^dF_2\\
 &\le&\dfrac{c}{\widehat{\mathbf{n}}h^d}\sum_{\mathbf{i}\ne\mathbf{j}}\left|K\left(\frac{\frac{\mathbf{i}}{\mathbf{n+1}}-\mathbf{s}_{0}}{h}\right)-K\left(\frac{\frac{\mathbf{j}}{\mathbf{n+1}}-\mathbf{s}_{0}}{h}\right)\right|K\left(\frac{\frac{\mathbf{j}}{\mathbf{n+1}}-\mathbf{s}_{0}}{h}\right)|Cov\left(\epsilon_{\mathbf{i}}^*,\epsilon_{\mathbf{j}}^*\right)|\\
  &&+\dfrac{c}{\widehat{\mathbf{n}}h^d}\sum_{\mathbf{i}\ne\mathbf{j}}\left[K\left(\frac{\frac{\mathbf{j}}{\mathbf{n+1}}-\mathbf{s}_{0}}{h}\right)\right]^2|Cov\left(\epsilon_{\mathbf{i}}^*,\epsilon_{\mathbf{j}}^*\right)|\\
 &\le&\dfrac{c_2}{\widehat{\mathbf{n}}h^d}\sum_{\mathbf{j}\in{\mathcal I}_{\mathbf{n}}}K\left(\frac{\frac{\mathbf{j}}{\mathbf{n+1}}-\mathbf{s}_{0}}{h}\right)\left(\sum_{\stackrel{\mathbf{i}\in{\mathcal I}_{\mathbf{n}}}{\|\mathbf{i-j}\|>Q}}\dfrac{\left\|\frac{\mathbf{i}}{\mathbf{n+1}}-\frac{\mathbf{j}}{\mathbf{n+1}}\right\|}{h}\left[\alpha_{1,\infty}(\|\mathbf{i-j}\|)\right]^{\frac{1}{2}}\right.\\
 &&\left. \ \ \ \ \ \ \ \ \ \ \ \ \ \ \ \ \ \ \ \ \ \ \ \ \ \ \ \ \ \ \ \ \ \ \ \ \ \ \ \ \ \ \ \ \ \ \ \ \ \ \ \ \ \ \ \ \ \ \ \ +\ \frac{Q^{2d}}{\min_{k=1,\cdots,d}(n_k)h}\right)\\
 &&+\dfrac{c_2}{\widehat{\mathbf{n}}h^d}\sum_{\mathbf{j}\in{\mathcal I}_{\mathbf{n}}}\left[K\left(\frac{\frac{\mathbf{j}}{\mathbf{n+1}}-\mathbf{s}_{0}}{h}\right)\right]^2\left(\sum_{\stackrel{\mathbf{i}\in{\mathcal I}_{\mathbf{n}}}{\|\mathbf{i-j}\|>Q}}\left[\alpha_{1,\infty}(\|\mathbf{i-j}\|)\right]^{\frac{1}{2}}+Q^d\right)
 \end{eqnarray*}
 and by replacing $Q$ by its value, we obtain
 \begin{eqnarray*}
&&\widehat{\mathbf{n}}h^dF_2 \\
  &\leq&\dfrac{c_2}{\widehat{\mathbf{n}}h^d\min_{k=1,\cdots,d}(n_k)h}\sum_{\mathbf{j}\in{\mathcal I}_{\mathbf{n}}}K\left(\frac{\frac{\mathbf{j}}{\mathbf{n+1}}-\mathbf{s}_{0}}{h}\right)\left(\sum^{+\infty}_{t=1}t^{-(\theta-2d)/2}+(\log \widehat{\mathbf{n}})^2\right)\\
 &&+\dfrac{c_2}{\widehat{\mathbf{n}}h^d}\sum_{\mathbf{j}\in{\mathcal I}_{\mathbf{n}}}\left[K\left(\frac{\frac{\mathbf{j}}{\mathbf{n+1}}-\mathbf{s}_{0}}{h}\right)\right]^2\left(\sum^{+\infty}_{t=1}t^{-(\theta-2d+2)/2}+\log \widehat{\mathbf{n}}\right),
 \end{eqnarray*}
 where $c_2$ is some positive constant. 
 Finally, we conclude that
 \begin{eqnarray*}
 \lim_{\mathbf{n}\to +\infty}\dfrac{\widehat{\mathbf{n}}h^d}{\log \widehat{\mathbf{n}}}\sup_{\mathbf{s}_0\in[0,1]^d}Var\left(\widehat{r}(\mathbf{s}_{0})\right)\leq c_2\int K^{2}(\mathbf{u})d\mathbf{u}<+\infty.
 \end{eqnarray*}


\end{document}